\title{ ~~\\ On the distribution of the order and index of
$g({\rm mod~}p)$ over residue classes II}
\author{Pieter Moree}
\documentclass[12pt]{article}
\usepackage{amssymb, latexsym, amsfonts}
\def\@ptsize{2}
\newtheorem{Thm}{Theorem}
\newtheorem{Conj}{Conjecture}
\newtheorem{Lem}{Lemma}
\newtheorem{Def}{Definition}
\newtheorem{Cor}{Corollary}

\newtheorem{Prop}{Proposition}
\newcommand{\qed}{\hfill $\Box$}

\begin{document}
\date{}
\maketitle
{\def\thefootnote{}
\footnote{\noindent 
Max-Planck-Institut f\"ur Mathematik, 
Vivatsgasse 7, 
D-53111 Bonn, Deutschland, E-mail:
moree@mpim-bonn.mpg.de\\}
{\def\thefootnote{}
\footnote{{\it Mathematics Subject Classification (2000)}.
11N37, 11N69, 11R45}}

\begin{abstract}
\noindent For a fixed rational number $g\not\in \{-1,0,1\}$ and
integers $a$ and $d$ we consider the
set  $N_g(a,d)$ of primes $p$ for which the order of $g({\rm mod~}p)$ is
congruent to $a({\rm mod~}d)$.
It is shown, assuming the
Generalized Riemann Hypothesis (GRH), that this set has a natural
density $\delta_g(a,d)$. Moreover, $\delta_g(a,d)$ is computed in terms of
degrees of certain Kummer extensions. Several properties of
$\delta_g(a,d)$ are established in case $d$ is a power of an odd prime.
\end{abstract}

\section{Introduction}
\label{introdu}
Let $g\not\in \{-1,0,1\}$ be a rational number (this 
assumption on $g$ will be maintained throughout this paper). For $u$
a rational number, let $v_p(u)$ denote the exponent of $p$ in the canonical
factorisation of $u$ (throughout the letter $p$ will be used to indicate
prime numbers). If $v_p(g)=0$, then there exists a smallest positive integer
$k$ such that $g^k\equiv 1({\rm mod~}p)$.
We put ord$_p(g)=k$.
This number is the {\it (residual)
order} of $g({\rm mod~}p)$.
The index of the subgroup generated by $g$ mod $p$ inside the multiplicative
group of residues mod $p$, $\left|({\mathbb Z}/p{\mathbb Z})^\times:
\langle g({\rm mod~}p)\rangle\right|$, is denoted by $r_g(p)$ and
called the {\it (residual) index} mod $p$ of $g$.
Note that
${\rm ord}_g(p)r_g(p)=p-1.$\\
\indent We let
$N_g(a_1,d_1;a_2,d_2)(x)$ count the number of  primes $p\le x$ with
$p\equiv a_1({\rm mod~}d_1)$ and ord$_g(p)\equiv a_2({\rm mod~}d_2)$.
For convenience denote $N_g(0,1;a,d)(x)$ by $N_g(a,d)(x)$. 
Although our main interest is in the behaviour of $N_g(a,d)(x)$ it turns out
that sometimes it is fruitful to partition $N_g(a,d)(x)$ in sets of the
form $N_g(a_1,d_1;a_2,d_2)(x)$ with a well-chosen $d_1$.
Let $r|s$ be positive integers. By $K_{s,r}$ we denote the number field $\mathbb Q(\zeta_{s},g^{1/r})$, 
where $\zeta_s=\exp(2\pi i/s)$.
The main result of this paper is as folows (for notational convenience $(a,b)$, $[a,b]$ will be written
for the greatest common divisor, respectively lowest common multiple of $a$ and $b$ and by (GRH) we
indicate that GRH is assumed).
\begin{Thm} 
\label{Main1}
{\rm (GRH)}.
Let $(a_1,d_1)=1$. Then
$$N_g(a_1,d_1;a_2,d_2)(x)=\delta_g(a_1,d_1;a_2,d_2){x\over \log x}+O_{d_1,d_2,g}({x\over 
\log^{3/2}x}),$$
for a number $\delta_g(a_1,d_1;a_2,d_2)$ that is given by {\rm (\ref{uitgebreid})}.
\end{Thm}
Specialising to $a_1=0$ and $d_1=1$ the following result is obtained.
\begin{Thm} {\rm (GRH)}.
\label{Main2}
We have
$$N_g(a,d)(x)=\delta_g(a,d){x\over \log x}+O_{d,g}({x\over 
\log^{3/2}x}),$$
with
\begin{equation}
\label{centraal}
\delta_g(a,d)=
\sum_{t=1\atop (1+ta,d)=1}^{\infty}
\sum_{n=1\atop (n,d)|a}^{\infty}{\mu(n)c_g(1+ta,dt,nt)\over 
[K_{[d,n]t,nt}:\mathbb Q]},
\end{equation}
where, for $(b,f)=1$,
$$c_g(b,f,v)=\cases{1 &if $\sigma_b|_{\mathbb Q(\zeta_{f})\cap K_{v,v}}={\rm id}$;\cr
0 & otherwise},$$
where $\sigma_b$ is the automorphism of $\mathbb Q(\zeta_f)$ that sends
$\zeta_f$ to $\zeta_f^b$.
\end{Thm}
In part I \cite{Moreealleen} this result was established, by a slightly different method, in case $3\le d\le 4$ and
explicit formulae for $\delta_g(a,3)$ and  $\delta_g(a,4)$ were derived. By a different method some 
of the
results for $d=4$ with weaker error term and restricted values of $g$ are established in
\cite{CM} and \cite{MC}, see also Corollary 1 of \cite{Moreealleen}.\\
\indent It turns out that the numbers $c_g$ appearing in (\ref{centraal}) have a strong
tendency to equal one. This motivates the following definition:
\begin{equation}
\label{nuldens}
\delta_g^{(0)}(a,d)=
\sum_{t=1\atop (1+ta,d)=1}^{\infty}
\sum_{n=1\atop (n,d)|a}^{\infty}{\mu(n)\over 
[K_{[d,n]t,nt}:\mathbb Q]}.
\end{equation}
For example, if $d|a$, then $\sigma_{1+ta}$ is the identity element of the Galois group of
$\mathbb Q(\zeta_{dt})$ over $\mathbb Q$ and then trivially $\delta_g(a,d)=\delta_g^{(0)}(a,d)$.\\
\indent Generically the degree $[K_{[d,n]t,nt}:\mathbb Q]$ 
appearing in (\ref{nuldens}) equals $\varphi([d,n]t)nt$. On
substituting this value in (\ref{nuldens}) a number $\delta(a,d)$ is obtained that no longer
depends on $g$:
\begin{equation}
\label{gemdens}
\delta(a,d)=
\sum_{t=1\atop (1+ta,d)=1}^{\infty}
\sum_{n=1\atop (n,d)|a}^{\infty}{\mu(n)\over 
\varphi([d,n]t)nt}.
\end{equation}
In \cite{Moreeaverage} it is shown that $\delta(a,d)$ is the average density of elements in a finite
field having order $\equiv a({\rm mod~}d)$. The number $\delta(a,d)$ can be a regarded as
a na\"{\i}ve heuristic for $\delta_g(a,d)$.\\
\indent In case $d$ is a power of an odd prime $q$ the coefficients $c_g$ are easily
evaluated. 
This case is considered in extenso in Section \ref{caseq}. 
There it is shown, for example, that $\delta_g(a,q^s)=\delta(a,q^s)$ for almost 
all integers $g$ with $|g|\le x$.\\

\section{Preliminaries}
\subsection{Some notation}
\label{notatie}
We recall some notation from our earlier papers \cite{Moreealleen, Moreeaverage}
on this topic.
For any Dirichlet character $\chi$ of $(\mathbb Z/d\mathbb Z)^*$, we let $h_{\chi}$ denote the
Dirichlet convolution of $\chi$ and $\mu$. 
As usual we let $L(s,\chi)$ denote the Dirichlet L-series for $\chi$.
For properties of $h_{\chi}$
the reader is referred to \cite{Moreealleen}. From \cite{Moreealleen} we furthermore
recall that
$$C_{\chi}(h,r,s)=
\sum_{(r,v)=1,~s|v}^{\infty}{h_{\chi}(v)(h,v)\over v\varphi(v)}
{\rm ~and~}A_{\chi}=\prod_{p:\chi(p)\ne 0}\left(1+{[\chi(p)-1]p\over [p^2-\chi(p)](p-1)}\right).$$
The constants $A_{\chi}$ turn out to be the basic constants for this
problem. 
In many cases $A_{\chi}\in \mathbb C\backslash \mathbb R$, see 
\cite[Table 3]{Moreeaverage}. It can be shown that
$$A_{\chi}\prod_{p|d}\left(1-{1\over p(p-1)}\right)=
A{L(2,\chi)L(3,\chi)\over L(6,\chi^2)}\prod_{r=1}^{\infty}\prod_{k=3r+1}^{\infty}
L(k,\chi^r)^{\lambda(k,r)},$$
where $A$ denotes the Artin constant and the numbers $\lambda(k,r)$ are non-zero integers that
can be related to Fibonacci numbers \cite{PFibo}. The latter expansion 
of $A_{\chi}$ can be used to approximate $A_{\chi}$ with high numerical accuracy 
(see \cite[Section 6]{Moreeaverage}).\\
\indent In \cite[Lemma 10]{Moreealleen} it is shown that $C_{\chi}(h,r,s)=cA_{\chi}$, where $c$ can
be explicitly written down and is in $\mathbb Q(\zeta_{o_{\chi}})$. Here $o_{\chi}$ is the
order of the character $\chi$, i.e. the smallest positive integer $k$ such that $\chi^k=\chi_0$, 
the trivial character.\\
\indent We recall from \cite{Moreealleen} that, if $(b,f)=1$,
\begin{equation}
\label{transitie}
\sum_{t\equiv b({\rm mod~}f)\atop t|v}\mu({v\over t})={1\over \varphi(f)}\sum_{\chi\in G_f}
{\overline{\chi(b)}}h_{\chi}(v),
\end{equation}
where the sum is over the characters in the character group $G_f$ of $(\mathbb Z/\mathbb Zf)^*$.

\subsection{Preliminaries on algebraic number theory}
We first review some properties of the Kronecker symbol (a rarely discussed
symbol in books on number theory). To this end we first recall the
definition of the Legendre and the Jacobi symbol. By definition the Legendre symbol $({n\over p})$, where
$p\ge 3$ is a prime number and $n\in \mathbb Z$, $p\nmid n$, is equal to $1$ if $n$
is a quadratic residue mod $p$, and to $-1$ otherwise.\\
\indent Let $m>0$ be an odd integer relatively prime to $n$. The Jacobi symbol $({n\over m})$ is defined as the product of the
Legendre symbols 
$({n\over m})=({n\over p_1})\cdots ({n\over p_s})$, where
$m=p_1\cdots p_s$ and each $p_i$ is a prime.\\
\indent The Kronecker symbol $({c\over d})$ is defined for $c\in \mathbb Z$, $c\equiv 0({\rm mod~}4)$ or
$c\equiv 1({\rm mod~}4)$, $c$ not a square, and $d\ge 1$ an integer; if $b=p_1p_2\cdots p_s$ is the
decomposition of $b$ as a product of primes, we put $({a\over -b})=({a\over b})=
({a\over p_1})({a\over p_2})\cdots ({a\over p_s})$. If $p$ is an odd prime $({a\over p})=0$ when $p$ divides
$a$, while $({a\over p})$ is the Legendre symbol $({a\over p})$ when $p$ does not divide $a$ and, finally,
$({a\over 2})=1$ when $a\equiv 1({\rm mod~}8)$, while  $({a\over 2})=-1$ when $a\equiv 5({\rm mod~}8)$.
Then if $a$ and $b$ are such that both the Jacobi and Kronecker symbols are defined, then these symbols
coincide. If $a$ is odd, then $({a\over 2})=$ Jacobi symbol $({2\over |a|})$. If $b>0$, $(a,b)=1$, $a$ is
odd, then $({a\over b})=({b\over |a|})$, where the symbol on the right hand side is the Jacobi symbol.
Most importantly, if $b>0$, $(a,b)=1$ and $a=2^r{\tilde a}$ with ${\tilde a}$ odd, then
$$({a\over b})=({2\over b})^r(-1)^{{{\tilde a}-1\over 2}{b-1\over 2}}({b\over |a|}),$$
where the symbols on the right hand side are Jacobi symbols.\\
\indent Let $K:\mathbb Q$ be an abelian number field. By the Kronecker-Weber theorem there exists an
integer $f$ such that $K\subseteq \mathbb Q(\zeta_f)$. The smallest such integer is called the
conductor of $K$. Note that $K\subseteq \mathbb Q(\zeta_n)$ iff $n$ is divisible by the conductor.
Note also that the conductor of a cyclotomic field is never congruent to $2({\rm mod~}4)$.
The following lemma allows one to determine all quadratic subfields of a given cyclotomic field
(for a proof see e.g. \cite[p. 263]{Weiss}).
\begin{Lem}
\label{conductor}
The conductor of a quadratic number field is equal to the absolute value of its discriminant.
\end{Lem} 
\indent Consider the cyclotomic extension $\mathbb Q(\zeta_f)$ of the rationals. There are $\varphi(f)$ distinct
automorphisms each determined uniquely by $\sigma_a(\zeta_f)=\zeta_f^a$, with $1\le a\le f$ and
$(a,f)=1$. We need to know when the restriction of such an automorphism to a given quadratic
subfield of $\mathbb Q(\zeta_f)$ is the identity. In this direction we have:
\begin{Lem}
\label{conductor2}
Let $\mathbb Q(\sqrt{d})\subseteq \mathbb Q(\zeta_f)$ be a quadratic field of discriminant $\Delta_d$ and
$b$ be an integer with $(b,f)=1$.
We have $\sigma_b|_{\mathbb Q(\sqrt{d})}={\rm id}$ iff $({\Delta_d\over b})=1$, with
$({\cdot \over \cdot})$ the Kronecker symbol.
\end{Lem}
{\it Proof}. Using Lemma \ref{conductor} we see that we can restrict to the case where $f=|\Delta_d|$.
Define $\chi$ by $\chi(b)=\sigma_b(\sqrt{d})/\sqrt{d}$, $1\le b\le |\Delta_d|$, $(b,\Delta_d)=1$. Then
$\chi$ is the unique non-trivial character of the character group of $\mathbb Q(\sqrt{d})$. As is
well-known (see e.g. \cite[p. 437]{Narkiewicz}), the primitive character induced by this is
$({\Delta_d\over b})$. Using Lemma \ref{conductor}, we see that $\chi$ is a primitive character mod
$|\Delta_d|$. Thus $\chi(b)=({\Delta_d\over b})$. Now
$\sigma_b|_{\mathbb Q(\sqrt{d})}={\rm id}$ iff $\chi(b)=({\Delta_d\over b})=1$. \qed\\

\noindent {\tt Remark 1}. Another way to prove Lemma \ref{conductor2} is to note that  
$\sigma_b|_{\mathbb Q(\sqrt{d})}={\rm id}$ iff there exists a prime $p\equiv b({\rm mod~}f)$ that
splits completely in $\mathbb Q(\sqrt{d})$. 
It is well-known that there exists a prime $p\equiv b({\rm mod~}f)$ that
splits completely in the field $\mathbb Q(\sqrt{d})$ 
iff $({\Delta_d\over p})=1$ (see e.g. \cite[p. 236]{Weiss}).
Since $({\Delta_d\over p})=({\Delta_d\over b})$, the result follows.\\

\noindent {\tt Remark 2}. The action of $\sigma$ on $\sqrt{d}$ can also be determined by
relating $\sqrt{d}$ to a Gauss sum `living' in $\mathbb Q(\zeta_f)$. It is straigthforward 
to determine the action of $\sigma$
on such a Gauss sum.

\subsection{Preliminaries on field degrees}
In order to explicitly evaluate certain densities in this paper, the 
following result will play a crucial role. 
Let $g_1\ne 0$ be a rational number. By $D(g_1)$ we denote the discriminant
of the field $\mathbb Q(\sqrt{g_1})$. 
The notation $D(g_1)$ along with the notation $g_0$,  
$h$ and $n_r$ introduced in the next lemma will reappear again and again in the sequel. 
\begin{Lem} {\rm \cite{Moreealleen}}.
\label{degree}
Write $g=\pm g_0^h$, where $g_0$ is
positive and not an
exact power of a rational. 
Let $D(g_0)$ denote the discriminant of the field $\mathbb Q(\sqrt{g_0})$.
Put $m=D(g_0)/2$ if $\nu_2(h)=0$ and $D(g_0)\equiv 4({\rm mod~}8)$
or $\nu_2(h)=1$ and $D(g_0)\equiv 0({\rm mod~}8)$, and
$m=[2^{\nu_2(h)+2},D(g_0)]$ otherwise. 
Put $$n_r=\cases{m &if $g<0$ and $r$ is odd;\cr
[2^{\nu_2(hr)+1},D(g_0)] &otherwise.}$$
We have
$$[K_{kr,k}:\mathbb Q]=[\mathbb Q(\zeta_{kr},g^{1/k}):\mathbb Q]={\varphi(kr)k\over
\epsilon(kr,k)(k,h)},$$
where, for $g>0$ or $g<0$ and $r$ even we have
$$\epsilon(kr,k)=\cases{2 &if $n_r|kr$;\cr
1 &if $n_r\nmid kr$,}$$
and for $g<0$ and $r$ odd we have
$$\epsilon(kr,k)=\cases{2 &if $n_r|kr$;\cr
 {1\over 2} &if $2|k$ and $2^{\nu_2(h)+1}\nmid k$;\cr
1 &otherwise.}$$
\end{Lem}
{\tt Remark}. Note that if $h$ is odd, then $n_r=[2^{\nu_2(r)+1},D(g)]$. Note that
$n_r=n_{\nu_2(r)}$.\\

\noindent From the latter lemma many consequences can be deduced.
\begin{Lem}
\label{graadoprekking}
Let $v,w$ and $z$ be natural numbers with $v|w$ and with $z$ an odd divisor of $w$.
Then
$[K_{zw,v}:\mathbb Q]=z[K_{w,v}:\mathbb Q]$.
\end{Lem}
{\it Proof}. The proof easily follows from Lemma \ref{degree} on observing that the odd part of
$n_r$ is squarefree and that $\varphi(zw)=z\varphi(w)$. \qed

\begin{Lem}
\label{opstart}
The intersection $\mathbb Q(\zeta_f)\cap K_{v,v}$ is equal to $\mathbb Q(\zeta_{(f,v)})$ or
a quadratic extension thereof.
More precisely,
$$[\mathbb Q(\zeta_f)\cap K_{v,v}:\mathbb Q(\zeta_{(f,v)})]={\epsilon([f,v],v)\over \epsilon(v,v)}.$$
\end{Lem}
{\it Proof}.
Clearly
this intersection field is abelian and contains $\mathbb Q(\zeta_{(f,v)})$. 
We have
\begin{equation}
\label{eerstedoorsnede}
[K_{[f,v],v}:\mathbb Q]={\varphi(f)[K_{v,v}:\mathbb Q]\over [\mathbb Q(\zeta_f)\cap K_{v,v}:\mathbb Q]}.
\end{equation}
On noting that $\varphi((f,v))\varphi([f,v])=\varphi(f)\varphi(v)$, it follows from Lemma \ref{degree} 
and (\ref{eerstedoorsnede}) that
\begin{equation}
\label{tweededoorsnede}
[\mathbb Q(\zeta_f)\cap K_{v,v}:\mathbb Q(\zeta_{(f,v)})]={[\mathbb Q(\zeta_f)\cap K_{v,v}:
\mathbb Q]\over \varphi((f,v))}=
{\epsilon([f,v],v)\over \epsilon(v,v)}.
\end{equation}
It is not difficult to infer from Lemma \ref{degree} that the latter quotient is either $1$ or $2$ (so
the apparent possibility $4$ does never arise). We conclude that 
$\mathbb Q(\zeta_f)\cap K_{v,v}$ is equal
to $\mathbb Q(\zeta_{(f,v)})$ or a quadratic extension thereof. \qed

\begin{Lem}
\label{degreeq}
Let $q$ be an odd prime, $s\ge 0$ and suppose $q\nmid v$. Put $q^*=({-1\over q})q$.
Consider the following conditions\\
i) $q|D(g_0)$, $s=0$ and ${n_1\over q}|v$.\\
ii) If $g<0$ and $2|v$, then $2^{\nu_2(h)+1}|v$.\\
We have
$$\mathbb Q(\zeta_{q^{s+1}})\cap K_{q^sv,q^sv}=\cases{\mathbb Q(\sqrt{q^*}) &if both i and ii are satisfied;\cr
\mathbb Q(\zeta_{q^s}) & otherwise.}$$
\end{Lem}
{\it Proof}. By Lemma \ref{opstart} we have $[\mathbb Q(\zeta_{q^{s+1}})\cap K_{q^sv,q^sv}:\mathbb 
Q(\zeta_{q^s})]=
\epsilon(q^{s+1}v,q^sv)/\epsilon(q^sv,q^sv)$
and hence, by Lemma \ref{degree}, $\mathbb Q(\zeta_{q^{s+1}})\cap K_{q^sv,q^sv}$ is a quadratic extenstion
of $\mathbb Q(\zeta_{q^0})=\mathbb Q$ if conditions i and ii are satisfied and
$\mathbb Q(\zeta_{q^{s+1}})\cap K_{q^sv,q^sv}=\mathbb Q(\zeta_{q^s})$ otherwise. Since 
$\mathbb Q(\sqrt{q^*})$ is the unique
quadratic subfield of $\mathbb Q(\zeta_q)$, the result now follows. \qed\\

\noindent {\tt Remark 3}. If $g>0$ or $g<0$ and $v$ is odd, then condition ii is vacuously satisfied.

\begin{Lem}
\label{dichtheidomschrijving}
Let $n$ be squarefree. Put $t_d=\prod_{p|(t,d)}p^{\nu_p(t)}$.
The density of primes $p$ such that $p\equiv 1+ta({\rm mod~}dt)$ and $p$ splits completely in $K_{nt,nt}$
equals zero if $(d,n)\nmid a$ or $(1+ta,d)>1$, otherwise it equals
\begin{equation}
\label{transfer}
{c_g(1+ta,dt,nt)\over [K_{[d,n]t,nt}:\mathbb Q]}
={c_g(1+ta,dt_d,nt)\over [K_{[d,n]t,nt}:\mathbb Q]}.
\end{equation}
\end{Lem}
{\it Proof}. This follows from Chebotarev's density
theorem together with the observation that the two systems of
congruences 
$$\cases{x\equiv 1+ta({\rm mod~}dt)\cr x\equiv 1({\rm mod~}nt)}{\rm ~and~}
\cases{x\equiv 1+ta({\rm mod~}dt_d)\cr x\equiv 1({\rm mod~}nt)}$$
are equivalent. \qed

\begin{Lem}
Assume that $(b,f_1f_2)=1$, $f_1|f_2$ and
\begin{equation}
\label{gelijk!}
{[K_{[f_1,v],v}:\mathbb Q]\over \varphi(f_1)}={[K_{[f_2,v],v}:\mathbb Q]\over \varphi(f_2)}.
\end{equation}
Then $c_g(b,f_1,v)=c_g(b,f_2,v)$.
\end{Lem}
{\it Proof}. By (\ref{eerstedoorsnede}) the assumption (\ref{gelijk!}) implies that
$$[\mathbb Q(\zeta_{f_1})\cap K_{v,v}:\mathbb Q]=[\mathbb Q(\zeta_{f_2})\cap K_{v,v}:\mathbb Q].$$ This,
together with the assumption that $f_1|f_2$ ensures that $\mathbb Q(\zeta_{f_1})\cap K_{v,v}=
\mathbb Q(\zeta_{f_2})\cap K_{v,v}=L$, say, whence $L=\mathbb Q(\zeta_{(f_1,f_2)})\cap K_{v,v}$. Since
the map $\sigma\in {\rm Gal}(\mathbb Q(\zeta_{f_1})/\mathbb Q)$ that sends $\zeta_{f_1}$ to $\zeta_{f_1}^b$
and the map  
$\sigma'\in {\rm Gal}(\mathbb Q(\zeta_{f_2})/\mathbb Q)$ 
that sends $\zeta_{f_2}$ to $\zeta_{f_2}^b$
act in the same way when restricted
to $\mathbb Q(\zeta_{(f_1,f_2)})$, it follows that $\sigma_b|_{\mathbb Q(\zeta_{f_1})\cap K_{v,v}}=
\sigma'_b|_{\mathbb Q(\zeta_{f_2})\cap K_{v,v}}$ and hence $c_g(b,f_1,v)=c_g(b,f_2,v)$. \qed

\subsection{Remaining preliminaries}
The following result is due to Wirsing \cite{Wirsing}.
\begin{Lem} {\rm \cite{Wirsing}}.
\label{wirsing}
Suppose $f(n)$ is a multiplicative function such that $f(n)\ge 0$, for
$n\ge 1$, and such that there are constants $\gamma_1$ and $\gamma_2$,
with $\gamma_2<2$, such that for every prime $p$ and every $\nu\ge 2$,
$f(p^{\nu})\le \gamma_1 \gamma_2^{\nu}$. Assume that as $x\rightarrow
\infty$,
$$\sum_{p\le x}f(p)\sim {\tau}{x\over \log x},$$
where $\tau>0$ is a constant. Then, as $x\rightarrow \infty$,
$$\sum_{n\le x}f(n)\sim {e^{-\gamma \tau}\over \Gamma(\tau)}
{x\over \log x}\prod_{p\le x}\left(1+{f(p)\over p}+{f(p^2)\over p^2}
+{f(p^3)\over p^3}+\cdots\right),$$
where $\gamma$ is Euler's constant and $\Gamma(\tau)$ denotes the
gamma-function.
\end{Lem}
We use it to establish the following lemma.
\begin{Lem}
\label{wirsinggevolg}
Let $d\ge 3$. The number of integers $1\le g\le x$ such that $D(g)$ has
no prime divisor $p$ with $p\equiv 1({\rm mod~}d)$ is
$O_d(x\log^{-1/\varphi(d)}x)$. The same assertion holds with
$D(g)$ replaced by $D(-g)$.
\end{Lem}
{\it Proof}. Denote the number of integers counted in the formulation
of the lemma by $T_d(x)$. We define a multiplicative function
$f_d(n)$ as follows:
$$
f_d(p^{\alpha})=\cases{0 &if $2\nmid \alpha$ and $p\equiv 1({\rm mod~}d)$;\cr
1 &otherwise.}
$$
Note that $T_d(x)=\sum_{g\le x}f_d(g)$.
Using Lemma \ref{wirsing}, it then follows that
$$T_d(x)=O_d\left({x\over \log x}\prod_{p\le x}\left(1+{f_d(p)\over p}
+{f_d(p^2)\over p^2}+\cdots \right)\right)
=O_d\left({x\over \log x}
\prod_{p\le x\atop p\not\equiv 1({\rm mod~}d)}(1+{1\over p})\right).$$
Since by Mertens' theorem for arithmetic progressions 
(cf. \cite{Williams}) we have
$$\prod_{p\le x\atop p\equiv a({\rm mod~}d)}(1+{1\over p})
\sim c_{a,d}\log^{1\over \varphi(d)}x,~x\rightarrow \infty,$$
for some $c_{a,d}>0$, the result follows. \qed

\section{The proofs of Theorems \ref{Main1} and \ref{Main2} }
In this section we analyse the growth behaviour of the counting 
functions $N_g(a,d)(x)$ and $N_g(a_1,d_1;a_2,d_2)(x)$. Throughout we restrict to those primes $p$ with
$\nu_p(g)=0$. Let $\omega(d)=\sum_{p|d}1$ denote the number of distinct prime divisors of $d$.\\
\indent Let  $$V_g(a,d;t)(x)=\#\{p\le x:r_g(p)=t,~p\equiv 1+ta({\rm mod~}dt)\}.$$
Note that $N_g(a,d)(x)=\sum_{t=1}^{\infty}V_g(a,d;t)(x)$.
If $(1+ta,d)>1$ then there is at most one prime counted by $V_g(a,d;t)(x)$ and this
prime has to divide $d$. It follows that 
$N_g(a,d)(x)=\sum_{t=1\atop (1+ta,d)=1}^{\infty}V_g(a,d;t)(x)+O(\omega(d))$.
Let $x_1=\sqrt{\log x}$. Assume GRH. By \cite[Lemma 7]{Moreealleen} it follows
that $=\#\{p\le x:r_g(p)>x_1\}=O_g(x\log^{-3/2}x)$. We thus infer that
\begin{equation}
\label{vanhieruit}
N_g(a,d)(x)=\sum_{t\le x_1\atop (1+ta,d)=1}V_g(a,d;t)(x)+O_{g,d}({x\over \log^{3/2}x}).
\end{equation}
For fixed $t$ the term $V_g(a,d;t)$ can be estimated by a variation of Hooley's classical
argument \cite{Hooley}. However, we need to carry this out with a certain uniformity. This requires
one to merely keep track of the dependence on $t$ of the various estimates. This results in
the following lemma.
\begin{Lem} {\rm (GRH)}.
\label{precies}
For $t\le x^{1/3}$ we have
$$V_g(a,d;t)(x)={x\over \log x}\sum_{n=1\atop (n,d)|a}^{\infty}{\mu(n)c_g(1+ta,dt,nt)\over 
[K_{[d,n]t,nt}:\mathbb Q]}+
O\left({x\log \log x\over \varphi(t)\log^2 x}
+{x\over \log^{2}x}\right),$$
where the implied constant depends at most on $d$ and $g$.
\end{Lem}
{\it Proof}. Let
$$M'_g(x,y)=\#\{p\le x:p\equiv 1+ta({\rm mod~}dt),~t|r_g(p),~qt\nmid r_g(p),~q\le y\}~{\rm ~and~}$$
$$M_g(x,y,z)=\#\{p\le x: qt|r_g(p),~y\le q\le z\},$$
where $q$ denotes a prime number.
Note that
$$V_g(a,d;t)(x)=M'_g(x,\xi_1)+
O(M_g(x,\xi_1,\xi_2))+O(M_g(x,\xi_2,\xi_3))
+O(M_g(x,\xi_3,{x-1\over t})).$$
We take $\xi_1={\log x/6}$, $\xi_2=\sqrt{x}\log^{-2}x$
and $\xi_3=\sqrt{x}\log x$. The three error terms were estimated in
\cite[Theorem 4]{Moreealleen}. Taking them together it is found that
\begin{equation}
\label{A}
V_g(a,d;t)(x)=M'_g(x,\xi_1)+O_{d,g}\left({x\log \log x\over \varphi(t)\log^2 x}
+{x\over \log^{2}x}\right),
\end{equation}
By inclusion and exclusion it follows that
$$
M'_g(x,\xi_1)=\sum_{P(n)\le \xi_1}\mu(n)\#\{p\le x:p\equiv 1+ta({\rm mod~}dt),~nt|r_g(p)\},
$$
where $P(n)$ denotes the greatest prime factor of $n$.
The integers $n$ counted in the latter sum are all less than $x^{1/3}$ (cf. (6) of
\cite{Hooley}). The counting functions in the latter sum can be estimated by 
an effective form of Chebotarev's
density theorem, cf. Theorem 3 of \cite{Moreealleen} and the discussion immediately
following this theorem. This yields that
$$\#\{p\le x:p\equiv 1+ta({\rm mod~}dt),~nt|r_g(p)\}={c_g(1+ta,dt,nt)\over [K_{[d,n]t,nt}:\mathbb Q]}
{\rm Li}(x)+O(\sqrt{x}\log x),$$
where the error term depends at most $d$ and $g$. Indeed, applying 
applying Theorem 3 of \cite{Moreealleen} results in an error term of
$O(\sqrt{x}\log(d_Lx^{[L:\mathbb Q]})/[L:\mathbb Q])$, with an absolute implied constant and
$L=K_{[d,n]t,nt}$. On invoking Lemma \ref{degree}
and \cite[Lemma 2]{Moreealleen} it follows that this is $O_{d,g}(\sqrt{x}\log x)$.
Proceeding as in section 6 of \cite{Hooley} it is then inferred that
\begin{equation}
\label{B}
M'_g(x,\xi_1)={x\over \log x}\sum_{n=1}^{\infty}{\mu(n)c_g(1+ta,dt,nt)\over 
[K_{[d,n]t,nt}:\mathbb Q]}+O_{d,g}({x\over \log^2 x}).
\end{equation}
If $nt|r_g(p)$, then $p\equiv 1({\rm mod~}nt)$. By the chinese remainder theorem it now follows
that if $(n,d)|a$ then $\#\{p\le x:p\equiv 1+ta({\rm mod~}dt),~nt|r_g(p)\}$ is finite for
every $x$ and hence $c_g(1+ta,dt,nt)=0$ (an alternative way to see this is
to note that in this case
$\sigma_{1+ta}$ does not act like the identity on $\mathbb Q(\zeta_{nt})$). It follows that in
the sum in (\ref{B}) we can restrict to those $n$ satisfying $(n,d)|a$.
On taking this into account and combining (\ref{A}) and (\ref{B}), the result follows. \qed\\

\noindent It is now straigthforward to establish Theorem \ref{Main2}.\\
{\it Proof of Theorem} \ref{Main2}.  
Recall that $x_1=\sqrt{\log x}$.
Combination of (\ref{vanhieruit}) and  Lemma \ref{precies}
yields
\begin{equation}
\label{C1}
N_g(a,d)(x)={x\over \log x}\sum_{t\le x_1\atop (1+ta,d)=1}\sum_{n=1\atop (n,d)|a}^{\infty}{\mu(n)c_g(1+ta,dt,nt)\over 
[K_{[d,n]t,nt}:\mathbb Q]}+O_{g,d}\left({x\over \log^{3/2}x}\right).
\end{equation}
Denote the latter double sum by $D(x)$. By Lemma \ref{degree} and \cite[Lemma 5]{Moreealleen} we find
$$D(x)=\delta_g(a,d)+O(\sum_{t>x_1}{h\over t\varphi(t)})=\delta_g(a,d)+O({h\over \sqrt{\log x}}).$$
On inserting the latter estimate in (\ref{C1}) the proof is then completed. \qed\\

\noindent A variation of the above (but notationally rather more awkward and hence we only sketch
it) gives Theorem \ref{Main1} with 
\begin{equation}
\label{uitgebreid}
\delta_g(a_1,d_1;a_2,d_2)=\sum_{t=1,~(1+ta_2,d_2)=1\atop 1+ta_2\equiv a_1({\rm mod~}(d_1,d_2t))}^{\infty}
\sum_{n=1}^{\infty}{\mu(n)c_g(a_1,d_1,1+ta_2,d_2t,nt)\over 
[K_{nt,nt}(\zeta_{d_1},\zeta_{d_2t}):\mathbb Q]},
\end{equation}
where, for $(b_1,f_1)=(b_2,f_2)=1$ and $b_1\equiv b_2({\rm mod~}(f_1,f_2))$, we define
$$c_g(b_1,f_1,b_2,f_2,v)=\cases{1 &if $\tau|_{\mathbb Q(\zeta_{[f_1,f_2]})\cap K_{v,v}}={\rm id}$;\cr
0 & otherwise},$$
where $\tau$ is the (unique) automorphism of $\mathbb Q(\zeta_{[f_1,f_2]})$ determined by
$\tau(\zeta_{f_1})=\zeta_{f_1}^{b_1}$ and $\tau(\zeta_{f_2})=\zeta_{f_2}^{b_2}$.\\

\noindent {\it Proof of Theorem} \ref{Main1}. This is a variation of the proof of Theorem \ref{Main2}. Most error terms can be
estimated as before on dropping the condition that $p\equiv a_1({\rm mod~}d_2)$, which
brings us in the situation of Theorem \ref{Main2}.\\ 
\indent We first generalise Lemma \ref{precies}.
For that we only have to replace $M'_g(x,\xi_1)$
 $M''_g(x,\xi_1)$ say, where $M''_g(x,\xi_1)$ is defined as $M'_g(x,\xi_1)$ with
$a=a_2$ and $d=d_2$, but where now
furthermore the primes $p$ are required to satisfy $p\equiv a_1({\rm mod~}d_1)$.
The estimation of $M''_g(x,\xi_1)$ can be carried out completely similarly to that
of $M'_g(x,\xi_1)$.\\
\indent The analogue of (\ref{vanhieruit}) is easily derived to be
$$
N_g(a_1,d_1;a_2,d_2)(x)=\sum_{t\le x_1,~(1+ta_2,d_2)=1\atop 1+ta_2\equiv a_1({\rm mod~}(d_1,d_2t))}
V_g(a_1,d_1;a_2,d_2;t)(x)+O_{g,d}({x\over \log^{3/2}x}).
$$
From here on the proof is completed as before. \qed

\section{The case where $d$ is an odd prime}
\label{caseq}
Let $d=q$ be an odd prime.
In this case it turns out to be fruitful to consider separately the
primes $p$ with $p\equiv 1({\rm mod~}q)$ and those with 
$p\not\equiv 1({\rm mod~}q)$.
\subsection{The case where $q|a$} 
\indent Trivially $N_g(a,q)=N_g(0,q)$ and w.l.o.g, we may assume
that $a=0$. Note that the primes counted by
$N_g(0,q)(x)$ must satisfy $p\equiv 1({\rm mod~}q)$. Let $j=\nu_q(p-1)$. Note
that $p\not\in N_g(0,q)$ iff $q^j|r_g(p)$. Thus we infer that
$$N_g(0,q)(x)=\#\{p\le x:p\equiv 1({\rm mod~}q)\}-$$
$$\sum_{j=1}^{\infty}\#\{p\le x:p\equiv 1({\rm mod~}q^j),~p\not\equiv 1({\rm mod~}q^{j+1}),~
q^j|r_g(p)\}.$$
The density of primes $p$ satisfying $p\equiv 1({\rm mod~}q^j)$, $p\not\equiv 1({\rm mod~}q^{j+1})$
and $q^j|r_g(p)$ can be computed by Chebotarev's density theorem and equals
$${1\over [K_{q^j,q^j}:\mathbb Q]}-{1\over [K_{q^{j+1},q^j}:\mathbb Q]}.$$
A more refined analysis \cite{Wiertelak1}, cf. \cite{Odoni, pappa2} (with weaker error term), shows that
$$N_g(0,q)(x)=\delta_g(0,q){\rm Li}(x)+O_{g,q}\left({x(\log \log x)^4\over \log^3 x}\right),$$
with
\begin{equation}
\label{classico}
\delta_g(0,q)={1\over q-1}-\sum_{j=1}^{\infty}
\left({1\over [K_{q^j,q^j}:\mathbb Q]}-{1\over [K_{q^{j+1},q^j}:\mathbb Q]}\right).
\end{equation}
Note that $N_g(1,q;0,q)(x)=N_g(0,q)(x)$ and hence $\delta_g(1,q;0,q)=\delta_g(0,q)$.
The density $\delta_g(0,q)$ can be explicitly evaluated using Lemma \ref{degree}: 
\begin{equation}
\label{gadoor}
\delta_g(1,q;0,q)=\delta_g(0,q)={q^{1-\nu_q(h)}\over q^2-1}.
\end{equation}
Let us now assume GRH. Using Theorem \ref{Main2} it is inferred that
$$\delta_g(0,q)=\delta_g^{(0)}(0,q)=\sum_{t=1}^{\infty}\sum_{n=1}^{\infty}
{\mu(n)\over [K_{[q,n]t,nt}:\mathbb Q]}=S_1+S_q,$$
say, where in $S_1$ we take together those $n$ with $q\nmid n$ and 
in $S_q$ those with $q|n$. We have
$$S_1=\sum_{v=1}^{\infty}{\sum_{n|v,~q\nmid n}\mu(n)\over [K_{qv,v}:\mathbb Q]}
=\sum_{j=0}^{\infty}{1\over [K_{q^{j+1},q^j}:\mathbb Q]},$$
and similarly
$$S_q=\sum_{t=1}^{\infty}\sum_{n=1\atop q|n}^{\infty}
{\mu(n)\over [K_{nt,nt}:\mathbb Q]}=-\sum_{t=1}^{\infty}\sum_{n=1\atop q\nmid n}^{\infty}
{\mu(n)\over [K_{qnt,qnt}:\mathbb Q]}=-\sum_{j=1}^{\infty}{1\over [K_{q^j,q^j}:\mathbb Q]},$$
On adding $S_q$ to $S_1$ we find (\ref{classico}) on noting that $[K_{q,1}:\mathbb Q]=q-1$.

\subsection{The case where $q\nmid a$}
In the remainder of this section we assume that $q\nmid a$.
Recall that $q^*=({-1\over q})q$.
\begin{Prop} {\rm (GRH)}.
\label{nulmodq}
Let $q$ be an odd prime and $q\nmid a$. Then
$$\delta_g(1,q;a,q)={1\over (q-1)^2}
\left(1-{q^{1-\nu_q(h)}\over q+1}\right).$$
In particular $\delta_g(1,q;a,q)\in \mathbb Q_{>0}$ and does not depend on $a$.
\end{Prop}
{\it Proof}. The density $\delta_g(1,q;a,q)$ equals, using (\ref{uitgebreid}) and Lemma \ref{dichtheidomschrijving},
$$\sum_{t=1,~q|t\atop (1+ta,q)=1}^{\infty}
\sum_{n=1\atop q\nmid n}^{\infty}{\mu(n)c_g(1+ta,qt,nt)\over
[K_{qnt,nt}:\mathbb Q]}=\sum_{t=1\atop q|t}^{\infty}
\sum_{n=1\atop q\nmid n}^{\infty}{\mu(n)c_g(1+ta,q^{1+\nu_q(t)},nt)\over
[K_{qnt,nt}:\mathbb Q]}.$$
Suppose that $q\nmid n$. By Lemma \ref{degreeq} it follows
that  $\mathbb Q(\zeta_{q^{1+\nu_q(t)}})\cap K_{nt,nt}=\mathbb Q(\zeta_{q^{\nu_q(t)}})$.
Since $1+ta\equiv 1({\rm mod~}q^{\nu_q(t)})$,
the automorphism $\sigma_{1+ta}$ in Theorem \ref{Main2} acts like the identity on the
latter field intersection
and hence $c_g(1+ta,q^{1+\nu_q(t)},nt)=1$. We thus infer that
\begin{equation}
\label{whatever}
\delta_g(1,q;a,q)=\sum_{t=1\atop q|t}^{\infty}
\sum_{n=1\atop q\nmid n}^{\infty}{\mu(n)\over
[K_{qnt,nt}:\mathbb Q]}.
\end{equation}
In particular, it follows that $\delta_g(1,q;a,q)$ does not depend on $a$. We present two ways 
to complete the proof from this point onwards.\\
\noindent {\it First way}. From (\ref{whatever}) we infer that
$$\delta_g(1,q;a,q)=\sum_{t=1}^{\infty}
\sum_{n=1,~q\nmid n}^{\infty}{\mu(n)\over
[K_{q^2nt,qnt}:\mathbb Q]}=\sum_{v=1}^{\infty}{\sum_{n|v,~q\nmid n}\mu(n)\over [K_{q^2v,v}:\mathbb Q]}
=\sum_{j=1}^{\infty}{1\over [K_{q^{j+1},q^j}:\mathbb Q]}.$$
Using Lemma \ref{degree} the latter sum is easily evaluated.\\
\noindent {\it Second way}. Note that $\sum_{0\le a\le q-1}\delta_g(1,q;a,q)$ equals the density
of primes $p$ with $p\equiv 1({\rm mod~}q)$ and hence
\begin{equation}
\label{sommetje}
\sum_{0\le a\le q-1}\delta_g(1,q;a,q)={1\over q-1}.
\end{equation}
Since, provided that $q\nmid a$, $\delta_g(1,q;a,q)$ does not depend on $a$,
we conclude from (\ref{sommetje}) that
$$\delta_g(1,q;a,q)={1\over q-1}\left({1\over q-1}-\delta_g(1,q;0,q)\right).$$
Now invoke (\ref{gadoor}). \qed\\

\noindent {\tt Remark 4}. Using a different method in \cite[Theorem 10]{Moreealleen} the values
of $\delta_g(1,3^s;a,3)$ for $s\ge 1$ were calculated.\\
 
\noindent In order to determine $\delta_g(a,q)$ it turns out to be convenient to determine
$\delta_g(a,d)-\delta_g(1,q;a,q)$ first.
\begin{Lem}
\label{q}
{\rm (GRH)}. Let $q$ be an odd prime and $q\nmid a$. Then
$$\delta_g(a,q)-\delta_g(1,q;a,q)=$$
$${1\over q-1}
-\sum_{v=1\atop q\nmid v}^{\infty}{\sum_{t\equiv -{1\over a}({\rm mod~}q),~t|v}\mu({v\over t})
\over [K_{qv,v}:\mathbb Q]}
-\sum_{v=1,~q\nmid v\atop \sqrt{q^*}\in K_{v,v}}^{\infty}{\sum_{({ta+1\over q})=-1,~t|v}\mu({v\over t})
\over [K_{qv,v}:\mathbb Q]}.$$
\end{Lem}
{\it Proof}. By Theorem \ref{Main1}, Theorem \ref{Main2} and Lemma \ref{dichtheidomschrijving}
we infer that
$$\delta_g(a,q)-\delta_g(1,q;a,q)=\sum_{t=1,~q\nmid t\atop (1+ta,q)=1}^{\infty}
\sum_{n=1\atop q\nmid n}^{\infty}{\mu(n)c_g(1+ta,q,nt)\over
[K_{qnt,nt}:\mathbb Q]}.$$
Let us restrict now to values of $n$ and $t$ that occur in the latter double sum.
We have, by Lemma \ref{opstart},
$$\mathbb Q(\zeta_q)\cap K_{nt,nt}=\cases{\mathbb Q(\sqrt{q^*}) & if $\sqrt{q^*}\in K_{nt,nt}$;\cr
\mathbb Q & otherwise.}$$
Using Lemma \ref{conductor2} it then follows that
$$c_g(1+ta,q,nt)=\cases{(1+({q^*\over 1+ta}))/2 & if $\sqrt{q^*}\in K_{nt,nt}$;\cr
1 & otherwise.}$$
By the properties of the
Kronecker symbol we have $(q^*/1+ta)=(1+ta/q)$, where the symbol $(1+ta/q)$
is just the Legendre symbol. We can thus write 
$\delta_g(a,q)-\delta_g(1,q;a,q)=J_1-J_2$, where
\begin{equation}
\label{negen}
J_1=\sum_{t=1,~q\nmid t\atop (1+ta,q)=1}^{\infty}\sum_{n=1\atop q\nmid n}^{\infty}
{\mu(n)\over [K_{qnt,nt}:\mathbb Q]}{\rm ~and~}
J_2=\sum_{t=1,~q\nmid t\atop ({1+ta\over q})=-1}^{\infty}\sum_{n=1,~q\nmid n\atop \sqrt{q^*}\in K_{nt,nt}}^{\infty}
{\mu(n)\over [K_{qnt,nt}:\mathbb Q]}.
\end{equation}
On writing $nt=v$ we obtain
\begin{eqnarray}
\label{tien}
J_1&=&\sum_{v=1\atop q\nmid v}^{\infty}{\sum_{t|v}\mu({v\over t})
\over [K_{qv,v}:\mathbb Q]}-\sum_{v=1\atop q\nmid v}^{\infty}{\sum_{t\equiv -{1\over a}({\rm mod~}q),~t|v}\mu({v\over t})
\over [K_{qv,v}:\mathbb Q]}\cr
&=&{1\over q-1}-\sum_{v=1\atop q\nmid v}^{\infty}{
\sum_{t\equiv -{1\over a}({\rm mod~}q),~t|v}\mu({v\over t})
\over [K_{qv,v}:\mathbb Q]}
\end{eqnarray}
and
$$J_2=\sum_{v=1,~q\nmid v\atop \sqrt{q^*}\in K_{v,v}}^{\infty}{\sum_{({1+ta\over q})=-1,~t|v}\mu({v\over t})
\over [K_{qv,v}:\mathbb Q]}.$$
On combining these expressions with $\delta_g(a,q)-\delta_g(1,q;a,q)=J_1-J_2$,
the result follows. \qed\\

\noindent {\tt Example}. For $q=3$ and $3\nmid a$ we obtain, on GRH, that
\begin{eqnarray}
\delta_g(2,3;1,3)&=&{1\over 2}
-\sum_{v=1\atop 3\nmid v}^{\infty}{\sum_{t\equiv 2({\rm mod~}3),~t|v}\mu({v\over t})
\over [K_{3v,v}:\mathbb Q]}
-\sum_{v=1,~3\nmid v\atop \sqrt{-3}\in K_{v,v}}^{\infty}{\sum_{t\equiv 1({\rm mod~}3),~t|v}\mu({v\over t})
\over [K_{3v,v}:\mathbb Q]}.\cr
&=&{1\over 2}-\sum_{v=1\atop 3\nmid v}^{\infty}{\sum_{t|v}\mu({v\over t})
\over [K_{3v,v}:\mathbb Q]}
+\sum_{v=1,~3\nmid v\atop \sqrt{-3}\not\in K_{v,v}}^{\infty}{\sum_{t\equiv 1({\rm mod~}3),~t|v}\mu({v\over t})
\over [K_{3v,v}:\mathbb Q]}.\cr
&=&\sum_{v=1,~3\nmid v\atop \sqrt{-3}\not\in K_{v,v}}^{\infty}{\sum_{t\equiv 1({\rm mod~}3),~t|v}\mu({v\over t})
\over [K_{3v,v}:\mathbb Q]}.\nonumber
\end{eqnarray}
More generally, we have 
$$\delta_g(2,3;a,3)=\sum_{v=1,~3\nmid v\atop \sqrt{-3}\not\in K_{v,v}}^{\infty}{\sum_{t\equiv a({\rm mod~}3),~t|v}\mu({v\over t})
\over [K_{3v,v}:\mathbb Q]}.$$
Rewriting this expression in terms of $h_{\chi}'s$ we obtain
Theorem 11 of \cite{Moreealleen}. \\ 

\noindent We now have the ingredients to establish the following result.
\begin{Thm} {\rm (GRH)}.
\label{expli}
Let $q$ be an odd prime and $q\nmid a$. Then
\begin{equation}
\label{C}
\delta_g^{(0)}(a,q)={1\over q-1}+{1\over (q-1)^2}
\left(1-{q^{1-\nu_q(h)}\over q+1}\right)
-\sum_{v=1\atop q\nmid v}^{\infty}{\sum_{t\equiv -{1\over a}({\rm mod~}q),~t|v}\mu({v\over t})
\over [K_{qv,v}:\mathbb Q]}
\end{equation}
and
\begin{equation}
\label{D}
\delta_g(a,q)=\delta_g^{(0)}(a,q)
-\sum_{v=1,~q\nmid v\atop \sqrt{q^*}\in K_{v,v}}^{\infty}{\sum_{({ta+1\over q})=-1,~t|v}\mu({v\over t})
\over [K_{qv,v}:\mathbb Q]}.
\end{equation}
If $q\nmid D(g_0)$, then $\delta_g(a,q)=\delta_g^{(0)}(a,q)$.
\end{Thm}
{\it Proof}. Considering the terms with $q|t$ and $q\nmid t$ in the double sum for $\delta_g(a,q)$
separately we have, using (\ref{whatever}) and (\ref{negen}) that 
$\delta_g^{(0)}(a,q)=\delta_q(1,q;a,q)+J_1$. On using Proposition \ref{nulmodq}, (\ref{tien}) and
Lemma \ref{q} the first two assertions are established.\\
\indent If $q\nmid D(g_0)$, then by Lemma \ref{degreeq} there is no integer $v$ such that 
$q\nmid v$ and $\sqrt{q^*}\in K_{v,v}$ and hence the double sum in (\ref{D}) equals zero.\qed\\

\noindent {\tt Remark 5}. The double sum (\ref{D})  can be rewritten, using Lemma \ref{degreeq}, as
$$\sum_{v=1,~q\nmid v}^{\infty}
\left({2\over [K_{qv,v}:\mathbb Q]}-{2\over (q-1)[K_{v,v}:\mathbb Q]}\right)
\sum_{({ta+1\over q})=-1,~t|v}\mu({v\over t}).$$

\noindent {\tt Remark 6}. Note that the proof of Theorem \ref{expli} makes essential use of the law of quadratic
reciprocity (this law enters in the proof of Lemma \ref{q}).\\

\noindent Let $G_q$ be the character group associated to $(\mathbb Z/\mathbb Zq)^*$. Since the latter
group is abelian it follows that
$G_q\cong (\mathbb Z/\mathbb Zq)^*$. The principal character mod $q$ will be denoted by $\chi_0$. Recall
that $o_{\chi}$ denotes the order of $\chi$ considered as an element of $G_{q}$.\\
\indent Using (\ref{transitie}) and the following lemma, $\delta_g(a,q)$ and $\delta_g^{(0)}(a,q)$ can be
expressed as simple linear combinations of the constants $C_{\chi}$ introduced in Section \ref{notatie}. Each such
constant can be explicitly evaluated and is of the form $cA_{\chi}$ with $c\in \mathbb Q(\zeta_{o_{\chi}})$.
This allows one to explicitly evaluate $\delta_g(a,q)$ and $\delta_g^{(0)}(a,q)$. For reasons of space
we only will work this out in the case $g\in {\cal G}$, where $\cal G$ is the set
of rational numbers $g$ that can not be written as $-g_0^h$ or $g_0^h$ with $h>1$ an integer and $g_0$ a
rational number.
\begin{Lem} {\rm \cite[Lemma 11]{Moreealleen}}.
\label{laatsteh}
Let $r,s$ be integers with $s|r$. Let $\chi$ be a 
Dirichlet character. Then, if $g>0$ or $g<0$ and $s$ is even,
$$\sum_{(r,v)=1}{h_{\chi}(v)\over [K_{sv,v}:\mathbb Q]}=
{1\over 
\varphi(s)}\left(C_{\chi}(h,r,1)+C_{\chi}(h,r,{n_s\over (n_s,s)})\right).$$
When $g<0$ and $s$ is odd, the latter sum equals
$$
{1\over 
\varphi(s)}\left(C_{\chi}(h,r,1)-{1\over 2}C_{\chi}(h,r,2)
+{1\over 2}C_{\chi}(h,r,2^{\nu_2(h)+1})+
C_{\chi}(h,r,{n_s\over (n_s,s)})\right).$$
\end{Lem}
Now we can formulate one of our main results.
\begin{Thm}
\label{bijnaa}
{\rm (GRH)}. Let $q$ be an odd prime and $q\nmid a$.
We have 
$$\delta_g(a,q)=\sum_{\chi\in G_q}\chi(-a)c_{\chi}A_{\chi},$$
where $c_{\chi}\in \mathbb Q(\zeta_{o_{\chi}})$ may depend on
$q$ and $g$ (but not $a$) and can be explicitly evaluated.
\end{Thm}
{\it Proof}. We only deal with the case where $g>0$ or $g<0$ and $h$ is odd (the remaining more space 
consuming case being left to the reader).\\ 
\indent Using the identity (\ref{transitie}) and Lemma \ref{laatsteh} we can rewrite (\ref{C}) as
$$
\delta_g^{(0)}(a,q)={1\over q-1}+$$ $${1\over (q-1)^2}
\left(1-{q^{1-\nu_q(h)}\over q+1}\right)
-{1\over (q-1)^2}\sum_{\chi\in G_q}\chi(-a)\left(C_{\chi}(h,q,1)+C_{\chi}(h,q,{n_q\over (n_q,q)})\right).$$
Similarly, using Remark 5, we rewrite (\ref{D}) as
$$\delta_g(a,q)=\delta_g^{(0)}(a,q)-{2\over (q-1)^2}\sum_{\chi\in G_q}\chi(-a)
\left(C_{\chi}(h,q,{n_q\over (n_q,q)})-C_{\chi}(h,q,n_1)\right)\sum_b {\overline{\chi(b)}},$$
where the sum is over the integers $1\le b\le q-1$ for which $({1-b\over q})=-1$.
Note that $n_1=n_q$. If $q\nmid D(g_0)$, the latter double sum equals zero and we infer (as before) that
$\delta_g(a,q)=\delta_g^{(0)}(a,q)$. If $q|D(g_0)$, then $C_{\chi}(h,q,n_1)=0$ and we infer that
$$\delta_g(a,q)=\delta_g^{(0)}(a,q)-{2\over (q-1)^2}\sum_{\chi\in G_q}\chi(-a)
C_{\chi}(h,q,{n_1\over q})\sum_b {\overline{\chi(b)}}.$$
By \cite[Lemma 10]{Moreealleen} we can write $C_{\chi}(h,r,s)$ as $cA_{\chi}$ with $c\in \mathbb Q(\zeta_{o_{\chi}})$, 
where $c$ can be explicitly given. Using this the proof is easily completed. \qed\\

\noindent The following result is an easy consequence of the latter proof.
\begin{Prop} {\rm (GRH)}.
If $h$ is odd and $8|D(g)$, then $\delta_g(a,q)=\delta_{-g}(a,q)$.
\end{Prop}
{\it Proof}. 
The assumptions imply that $|n_1|=|n_q|=|D(\pm g)|$ and on noting that $C_{\chi}(h,r,s)=C_{\chi}(h,r,-s)$ the
character sum expression for $\delta_g(a,q)$ given in the proof of Theorem \ref{bijnaa} is seen to equal that
of $\delta_{-g}(a,q)$ in case $q\nmid a$. If $q|a$, then by (\ref{gadoor}) we have
$\delta_{\pm g}(0,q)=q/(q^2-1)$. \qed\\

\noindent The following result demonstates Theorem \ref{bijnaa} in the special (but 
important) case where $g\in {\cal G}$. Since almost all integers are in $\cal G$ it can be thought of as the
set of `generic' integers $g$.
\begin{Prop}
\label{explicietq}
{\rm (GRH)}. Let $q$ be an odd prime and $q\nmid a$.
Suppose that $g\in G$.
Put $$\epsilon_g(\chi)=\cases{1 & if $2\nmid D(g)$;\cr
{\chi(2)\over 4} &if $4||D(g)$;\cr
{\chi(2)^2\over 16} & if $8|D(g)$.}$$
If $q\nmid D(g)$, then
$\delta_g(a,q)$ equals 
$\delta_q^0(a,q)$  which on its turn equals
$${q^2\over (q-1)(q^2-1)}-{1\over (q-1)^2}\sum_{\chi\in G_q}\chi(-a)A_{\chi}
\left(1+\epsilon_g(\chi)\prod_{p|2D(g)}{p(\chi(p)-1)\over  p^3-p^2-p+\chi(p)}\right).$$
If $q|D(g)$, then
$$\delta_g(a,q)={q^2\over (q-1)(q^2-1)}$$
$$-{1\over (q-1)^2}\sum_{\chi\in G_q}\chi(-a)A_{\chi}
\left(1+\epsilon_g(\chi)\Big\{1+2\sum_b {\overline{\chi(b)}}\Big\}\prod_{p|2D(g)/q}{p(\chi(p)-1)\over  
p^3-p^2-p+\chi(p)}\right),$$
where the sum is over all $1\le b\le q-1$ for which $({1-b\over q})=-1$.
\end{Prop}
{\it Proof}. By Lemma \ref{degree} we have
$$n_1=n_q=\cases{[2,D(g_0)] & if $g>0$;\cr
D(g_0)/2 & if $g<0$ and $D(g)\equiv 4({\rm mod~}8)$;\cr
[4,D(g_0)] & if $g<0$ and $D(g)\not\equiv 4({\rm mod~}8)$.}
$$
Now note that $n_1=n_q=[2,D(g)]$ and that $q|D(g_0)$ iff $q|D(g)$. Working out the formulae involving the $C_{\chi}$'s in
the proof of Theorem \ref{bijnaa} using \cite[Lemma 10]{Moreealleen}, the proof is
then completed.\qed\\

\noindent In case $d=q$ is an odd prime it can be shown \cite[Proposition 7]{Moreeaverage}
that
$$\delta(a,q)={q^2\over (q-1)(q^2-1)}-{1\over (q-1)^2}\sum_{\chi\in G_q}\chi(-a)A_{\chi}.$$
Assume GRH. Proposition \ref{explicietq}  shows that if $g\in G$ and $|g|$ tends to infinity, then
$\delta_g(a,q)$ tends to $\delta(a,q)$. Proposition \ref{explicietq} also shows that if
$g\in G$ and $D(g)$ contains a prime divisor $p$ with $p\equiv 1({\rm mod~}q)$, then
$\delta_g(a,q)=\delta(a,q)$. It is not difficult to show that almost all integers $g$ with
$|g|\le x$ satisfy $g\in G$ and have $D(g)$ with at least one prime divisor $p$
with $p\equiv 1({\rm mod~}q)$. 
Consequently we infer that for almost all integers $g$ with
$|g|\le x$, $\delta_g(a,q)=\delta(a,q)$. 
This shows that for a `generic' $g$, $\delta_g(a,q)=\delta(a,q)$. Furthermore, if
$\delta_g(a,q)$ is not equal to $\delta(a,q)$ then usually it will be quite close to it.
The results in the next section allow one to also draw these conclusions in the case 
where $d$
is an odd prime power, which will be done in the final section.

\section{The case where $d$ is an odd prime power}
The case where $d=q^s$ with $q$ an odd prime is easily reduced to the case $d=q$ by the following
result. 
\begin{Thm}
\label{locald}
{\rm (GRH)}. Suppose that $d|d_1$, the quotient $d_1/d$ is odd and $\omega(d_1)=\omega(d)$.
Then $\delta_g(a,d_1)={d\over d_1}\delta_g(a,d)$.
\end{Thm}
\begin{Cor} 
Let $q$ be an odd prime and $j\ge 1$, then $\delta_g(a,q^j)=q^{1-j}\delta_g(a,q)$.
\end{Cor}
\noindent {\tt Remark 7}. From formula (\ref{gemdens}) 
for $\delta(a,d)$ it is easily inferred that if $d|d_1$ and $\omega(d_1)=\omega(d)$, then
\begin{equation}
\label{delta}
 \delta(a,d_1)={d\over d_1}\delta(a,d).
\end{equation}
 
\noindent {\it Proof of Theorem} \ref{locald}. If $d|d_1$ and $\omega(d)=\omega(d_1)$ and for all
$n$ and $t$ with $(1+ta,d)=1$, $(n,d)|a$ and $n$ is squarefree, we have
$[K_{[d_1,n]t,nt}:\mathbb Q]/\varphi(d_1)=[K_{[d,n]t,nt}:\mathbb Q]/\varphi(d)$, then using 
Lemma \ref{gelijk!} we infer that
\begin{eqnarray}
\delta_g(a,d_1)&=&\sum_{t=1\atop (1+ta,d)=1}^{\infty}
\sum_{n=1\atop (n,d)|a}^{\infty}{\mu(n)c_g(1+ta,d_1t,nt)\over [K_{[d_1,n]t,nt}:\mathbb Q]}\nonumber\\
&=&{d\over d_1}\sum_{t=1\atop (1+ta,d)=1}^{\infty}
\sum_{n=1\atop (n,d)|a}^{\infty}{\mu(n)c_g(1+ta,dt,nt)\over [K_{[d,n]t,nt}:\mathbb Q]}={d\over d_1}
\delta_g(a,d),\nonumber
\end{eqnarray}
where we used that $\varphi(d_1)/\varphi(d)=d_1/d$. Now invoke Lemma \ref{graadoprekking}. \qed

\section{Connection between $\delta_g(a,q^s)$ and $\delta(a,q^s)$}
Define ${\overline \delta}_g(d)=(\delta_g(0,d),\delta_g(1,d),\dots,\delta_g(d-1,d))$ (if this exists)
and $${\overline \delta}(d)=(\delta(0,d),\delta(1,d),\dots,\delta(d-1,d)).$$
The next result implies that, under
GRH, for almost all integers $g$ we have ${\overline \delta}_g(q^s)={\overline \delta}(q^s)$.
\begin{Thm} {\rm (GRH)}.
\label{sseven}
Let $s\ge 1$ and $q$ an odd prime. Then there are at most $O_q(x\log^{-1/(q-1)}x)$ integers $g$
with $|g|\le x$ such that ${\overline \delta}_g(q^s)\ne {\overline \delta}(q^s)$.
\end{Thm}
{\it Proof}. By the results of the previous section it is enough to prove this in case $s=1$. 
By \cite[Proposition 7]{Moreeaverage} we have, for $1\le a\le q-1$,
$$\delta(a,q)={q^2\over (q-1)(q^2-1)}-{1\over (q-1)^2}\sum_{\chi\in G_q}\chi(-a)A_{\chi},$$
and $\delta(0,q)=q/(q^2-1)$ \cite[Theorem 1]{Moreeaverage}.
Since there
are at most $O(\sqrt{x}\log x)$ integers $|g|\le x$ that are not in $\cal G$, we can restrict ourselves to
the case where $g\in {\cal G}$. For such a $g$ we then have, using (\ref{gadoor}), that
$\delta_g(0,q)=\delta(0,q)$. Now by Proposition \ref{explicietq} we infer that if $D(g)$ has a prime divisor
$p$ with $p\equiv 1({\rm mod~}q)$, then ${\overline \delta}_g(q)={\overline \delta}(q)$ (since $\chi(p)=1$ for
every $\chi\in G_q$). It follows that the number of $g\in {\cal G}$ with $|g|\le x$ such that  
${\overline \delta}_g(q)\ne {\overline \delta}(q)$ is bounded above by the number of $g$ with $|g|\le x$ such
that $D(g)$ has no prime factor $p$ satisfying $p\equiv 1({\rm mod~}q)$.  
By Lemma \ref{wirsinggevolg} the proof is then completed. \qed\\

\noindent If it is not true that ${\overline \delta}_g(q^s)={\overline \delta}(q^s)$, then our final result
shows that ${\overline \delta}_g(q^s)$ will be close to ${\overline \delta}(q^s)$,
\begin{Prop} {\rm (GRH)}.
Suppose that $g\in {\cal G}$. As $|g|$ tends to infinity, ${\overline \delta}_g(q^s)$ tends to 
${\overline \delta}(q^s)$.
\end{Prop}
{\it Proof}. A simple consequence of the results in the previous section, 
the formula given for $\delta(a,q)$ in the proof of Theorem \ref{sseven}, (\ref{gadoor}) and
Proposition \ref{explicietq}. \qed\\

\noindent For a numerical illustration of the latter result (with $g=-19$) see Table 1 of
\cite{Moreeaverage}.

\end{document}